\documentclass[a4paper,12pt]{extarticle}
\usepackage{setspace}
\usepackage{mathtext}
\usepackage{graphicx}
\usepackage{amscd}
\usepackage[T2A]{fontenc}
\usepackage[english]{babel}
\usepackage{latexsym,amsfonts,amssymb,amsmath,longtable,amsthm}
\usepackage[pic]{xy}
\usepackage{bbm,makeidx}
\usepackage[centerlast,small]{caption}
\usepackage{amsmath}
\usepackage[cp1251]{inputenc}
\usepackage[usenames, dvipsnames]{color}
\newtheorem{lem}{{\scshape Lemma}}
\newtheorem{cor}{{\scshape Corollary}}
\newtheorem{rmk}{{\scshape Remark}}

\newtheorem{ttt}{{\scshape Theorem}}
\newtheorem{prp}{{\scshape Proposition}}

\newtheorem{problem}{{\scshape Problem}}
\newtheorem{con}{{\scshape Conjecture}}

\begin{document}
\title{On groups where the twisted conjugacy class of the unit element is a subgroup}
\author{D.~L.~Gon\c{c}alves\footnote{The author is partially supported by FAPESP-Funda\c c\~ao 
de Amparo a Pesquisa do Estado de S\~ao Paulo, Projeto Tem\'atico 
Topologia Alg\'ebrica, Geom\'etrica 2012/24454-8.}, T.~Nasybullov\footnote{The author is supported by the Research Foundation -- Flanders (FWO): postdoctoral grant  12G0317N and travel grant  V417817N.}}
\date{}
\maketitle
\begin{abstract}
We study groups $G$ where the $\varphi$-conjugacy class $[e]_{\varphi}=\{g^{-1}\varphi(g)~|~g\in G\}$ of the unit element is a subgroup of $G$ for every automorphism $\varphi$ of $G$. If $G$ has $n$ generators, then we prove that the $k$-th member of the lower central series has a finite verbal width bounded in terms of $n,k$. Moreover, we prove that if such group $G$ satisfies the descending chain condition for normal subgroups, then $G$ is nilpotent, what generalizes the result from \cite{BarNasNes}. Finally, if $G$ is a finite abelian-by-cyclic group, we construct a good upper bound of the nilpotency class of $G$.
~\\

\end{abstract}
\section{Introduction}
Conjugacy classes in a group reflect properties of this group. For example, the groups where every conjugacy class is finite (so called \textit{FC-groups} introduced by Baer in \cite{Bae})  have very restrictive structure. In particular, B.~Neumann discovered that if $G$ is a group in which all conjugacy classes are finite with bounded size, then the derived subgroup $G^{\prime}$ is finite \cite{Neu}. Another interesting connections between properties of conjugacy classes and properties of groups can be found, for example, in \cite{CamCam, Kel}.

The notion of $\varphi$-conjugacy classes generalizes the notion of conjugacy classes. Let $G$ be a group and $\varphi$ be an automorphism of $G$. Elements $x, y$ of the group $G$ are said to be \textit{(twisted) $\varphi$-conjugated} if there exists an element $z$ in $G$ such that $x=zy\varphi(z)^{-1}$. In particular, if $\varphi=id$, then we have the usual definition of conjugated elements. The relation of $\varphi$-conjugation is an equivalence relation, equivalence classes of which are called \textit{$\varphi$-conjugacy classes}. Denote by $[x]_{\varphi}$ the $\varphi$-conjugacy class of the element $x$ and by $R(\varphi)$ the number of these classes. The $\varphi$-conjugacy classes are also called \textit{Reidemeister classes}, and the number $R(\varphi)$ is called \textit{the Reidemeister number} of the automorphism $\varphi$. As well as conjugacy classes, twisted conjugacy classes in a group reflect properties of the group itself. For example, the following conjecture (partially proved in \cite[Proposition 5.2]{FelNas} and \cite[Theorems A, B and C]{Jab}) describes very strong connection between the Reidemeister number $R(\varphi)$ and the structure of the finitely generated residually finite group. 

~\\
\noindent\textbf{\cite[Conjecture R]{FelTro}.} Let $G$ be a finitely generated residually finite group. If there exists an automorphism $\varphi$ such that $R(\varphi)$ is finite, then $G$ is almost solvable.

~\\
Another interesting conjecture states the connection between the twisted conjugacy class $[e]_{\varphi}$ of the unit element $e$ and the structure of the group itself.

~\\
\noindent\textbf{\cite[Conjecture 1]{BarNasNes}.} Let $G$ be a group such that the class $[e]_{\varphi}$ is a subgroup of $G$ for every $\varphi\in {\rm Aut}~G$. Then $G$ is nilpotent. 

~\\
In \cite{BarNasNes} it is proved that if the class $[e]_{\varphi}$ is a subgroup of $G$ for every inner automorpism $\varphi$, then $G$ belongs to the Kurosh-Chernikov class $\overline{Z}$ and there exists a strictly descending series of normal subgroups in $G$. In particular, if $G$ is a finite group in which the class $[e]_{\varphi}$ is a subgroup for every inner automorphism $\varphi$, then $G$ is nilpotent.

In the present paper we continue to study groups where the twisted conjugacy class $[e]_{\varphi}$ is a subgroup of $G$ for every automorphism $\varphi$. In particular, we study the lower central series of such groups and prove that if $G$ is a group with $n$ generators such that the class $[e]_{\varphi}$ is a subgroup of $G$ for every inner automorphism $\varphi$, then the $k$-th member of the lower central series of $G$ can be decomposed into the product of the finite number (depending on $n$ and $k$) of subgroups of the form $[e]_{\varphi_i}$ for different inner automorphisms $\varphi_i$ (Theorem \ref{wfin}). 
This result has several interesting corollaries.

The first important corollary of Theorem \ref{wfin} is about the verbal width of the members of the lower central series.
Let $F_n$ be a free group with generators $x_1$, $\dots$, $x_n$ and $w(x_1,\dots,x_n)$ be an element of $F_n$. For the elements $g_1,\dots, g_n$ of a group $G$ we use the symbol $w(g_1,\dots,g_n)$ to denote the value of $w$ on the elements $g_1,\dots,g_n$. Denote by $w(G)=\langle w(g_1,\dots,g_n)~|~g_1,\dots,g_n\in G\rangle$ the subgroup of $G$ defined by the word $w\in F_n$. For the element $g\in w(G)$ denote by $l_w(g)$ the minimal number $k$ such that $g$ is a product of $k$ elements of the form $w(g_1,\dots,g_n)^{\pm1}$ for $g_1,\dots,g_n\in G$. The supremum of the values $l_w(g)$ for $g\in w(G)$ is called
\textit{the (verbal) width of $w(G)$} and is denote by ${\rm wid}(w(G))$. Width of verbal subgroups in different groups has been studied by various authors (see, for example, \cite{Seg} and and references therein). It is well known that
the commutator width (the word $w$ is a commutator) of a free non-abelian group is infinite  \cite{Bar, Rhe}, but the commutator width of a finitely generated nilpotent group is finite \cite{AllRom}. In the present paper we prove that if $G$ is a group with $n$ generators such that the class $[e]_{\varphi}$ is a subgroup of $G$ for every inner automorphism $\varphi$, then the width of the $k$-th member of the lower central series $\gamma_k(G)$ is finite and is bounded in terms of $n,k$.

As an another interesting corollary of Theorem \ref{wfin}, we prove that if $G$ is a finitely generated group which satisfies the descending chain condition for normal subgroups such that $[e]_{\varphi}$ is a subgroup of $G$ for every inner automorphism $\varphi$, then $G$ is nilpotent (Corollary \ref{minn}). This result generalises the result from \cite[Theorem 3]{BarNasNes}. 

As we already said, if $G$ is a finite group such that the class $[e]_{\varphi}$ is a subgroup of $G$ for every $\varphi$, then $G$ is nilpotent \cite[Theorem 3]{BarNasNes}. In the present paper for finite groups of the form $1\to A\to G\to \mathbb{Z}_n\to1$, where $A$ is abelian and $n$ is odd, we construct a good upper bound of the nilpotency class of $G$ which depends only on the number $n$ (Theorem \ref{corex}). 

Several open problems are formulated throughout the paper.

The authors are greatful to Vasiliy Bludov and Evgeny Khukhro for the examples from Section \ref{exa}.
\section{Definitions and known results}\label{prem}
We use the following notation. If $x,y$ are element of a group $G$, then $x^y=y^{-1}xy$ is \textit{conjugacy of $x$ by $y$} and $[x,y]=x^{-1}y^{-1}xy$ is \textit{a commutator of $x$ and $y$}. We use the symbol $\gamma_k(G)$ to denote the \textit{$k$-th term of the lower central series}: $\gamma_1(G)=G$, $\gamma_{k+1}(G)=[\gamma_k(G),G]$. The group $G$ is called \textit{nilpotent group of nilpotency class $m$} if $\gamma_{m+1}(G)=e$ and $\gamma_k(G)\neq e$ for all $k\leq m$. The group $G$ is called \textit{residually nilpotent} if $\bigcap_{k=1}^{\infty}\gamma_k(G)=e$. Every nilpotent group is obviously residually nilpotent.

Let $F_n$ be a free group with generators $x_1$, $\dots$, $x_n$ and $w$ be an element of $F_n$. If $G$ is a group and $g_1$, $\dots$, $g_n$ are elements of $G$, then we use the symbol $w(g_1,\dots,g_n)$ to denote the value of $w$ on the elements $g_1,\dots,g_n$. So we have a map $w:G^n\to G$ defined by $w$. The subgroup $w(G)=\langle w(g_1,\dots,g_n)~|~g_1,\dots,g_n\in G\rangle$ generated by all values of $w$ on $G$ is called \textit{a verbal subgroup of $G$ defined by the word $w\in F_n$}. The subgroup $\gamma_k(G)$ is a verbal subgroup defined by the word $\gamma_k$ from $F_k$: $\gamma_1(x_1)=x_1$, $\gamma_{k+1}(x_1,\dots,x_{k+1})=[\gamma_k(x_1,\dots,x_k),x_{k+1}]$. \textit{The (verbal) width ${\rm wid}(w(G))$} of the verbal subgroup $w(G)$ defined by the word $w$ is the smalles number $m\in\mathbb{N}\cup\{\infty\}$ such that every element of $w(G)$ is the product of at most $m$ words $w(g_1,\dots,g_n)$ and their inverses. The equality $w(G)=v(G)$ for different words $v,w$ in general does not imply the equality ${\rm wid}(w(G))={\rm wid}(v(G))$.

Let $\varphi$ be an automorphism of the group $G$. Elements $x,y$ of $G$ are said to be \textit{$\varphi$-conjugated} if there exists an element $z$ in $G$ such that $x=z^{-1}y\varphi(z)$. The relation of $\varphi$-conjugation is an equivalence relation, equivalence classes of which are called \textit{$\varphi$-conjugacy classes}. Denote by $[x]_{\varphi}$ the $\varphi$-conjugacy class of the element $x$. If $\hat{g}:x\mapsto x^g$ is an inner automorphism of $G$ induced by the element $g$, then for the sake of simplicity we will write $[x]_g$ instead of $[x]_{\hat{g}}$ in order to denote the $\hat{g}$-conjugacy class of the element $x$. In particular, the $\hat{g}$-conjugacy class of the unit element $e$ has the following form.
$$[e]_g=\{[x,g]~|~x\in G\}$$

The following two lemmas about the twisted conjugacy class $[e]_g$ are obvious.
\begin{lem}\label{easy1} If $G$ is a group such that $[e]_g$ is a subgroup of $G$ for every $g\in G$, then the element $e\neq g\in G$ does not belong to $[e]_g$.
\end{lem}
\begin{lem}\label{easy2} Let $G$ be a group such that $[e]_g$ is a subgroup of $G$ for every $g\in G$ and $H$ be a normal subgroup of $G$. If $g\in H$, then $[e]_g<H$.
\end{lem}

The following two lemmas are proved in \cite[Propositions 1, 3, 5]{BarNasNes}.
\begin{lem}\label{nsub} Let $G$ be a group and $\varphi$ be an automorphism of $G$. If the class $[e]_{\varphi}$ is a subgroup of $G$, then this subgroup is normal in $G$.
\end{lem}
\begin{lem}\label{prim} Let $\varphi$, $\psi$ be automorphisms of a group $G$ such that $[e]_{\varphi}$, $[e]_{\psi}$ are subgroups of $G$. Then $[e]_{\varphi\psi}\subseteq[e]_{\varphi}[e]_{\psi}$.
\end{lem}
\begin{lem}\label{quu} Let $G$ be a group and $[e]_g$ be a subgroup of $G$. Then for every normal subgroup $N$ of $G$ the class $[e]_{\overline{g}}$ is a subgroup of $G/N$, where $\overline{g}=gN\in G/N$.
\end{lem}
The word $w$ from $F_n=\langle x_1,\dots,x_n\rangle$ is called \textit{an outer commutator word} if it is obtained from a finite sequence of distinct variables $x_1,\dots,x_n$ using only commutator brackets. For example, the word $[[x_1,x_3],x_2]$ is an outer commutator word, while the words $[x_1,x_2][x_3,x_4]$ and $[[x_1,x_2],[x_1,x_3]]$ are not outer commutator words. The following lemma is proved in \cite[Lemma 2(1)]{BarNasNes}.
\begin{lem}\label{vnesh}Let $G$ be a group such that the twisted conjugacy class $[e]_{g}$ is a subgroup of $G$ for every $g\in G$. If $w\in F_n$ is an outer commutator word, then $w(g_1,\dots,g_n)$ belongs to $[e]_{g_1}\cap\dots\cap[e]_{g_n}$.  
\end{lem}
From Lemma \ref{easy1}, Lemma \ref{easy2} and Lemma \ref{vnesh} follows that 
$$[e]_{g_1}>[e]_{[g_1,g_2]}>[e]_{[g_1,g_2,g_3]}>\dots$$
 is a strictly descending chain of normal subgroups if each member of this chain is not trivial.
The following conjecture of Bardakov-Nasybullov-Neshchadim is formulated in \cite[Conjecture 1]{BarNasNes} and \cite[Problem 18.14]{Kou}.
\begin{con}\label{prob1}If $G$ is a group such that the class $[e]_{\varphi}$ is a subgroup of $G$ for every automorphism $\varphi$, then $G$ is nilpotent.
\end{con}
We say that the group $G$ satisfies \textit{the descending chain condition for normal subgroups (Min-n condition)} if for every chain $H_1>H_2>H_3>\dots$ of normal subgroups there exists a number $n$ such that $H_n=H_{n+1}=\dots$ We say that $G$ satisfies \textit{the ascending chain condition for normal subgroups (Max-n condition)} if for every chain $H_1<H_2<H_3<\dots$ of normal subgroups there exists a number $n$ such that $H_n=H_{n+1}=\dots$ The following result from \cite[Theorem 3]{BarNasNes} gives a partial positive solution of Conjecture \ref{prob1}.
\begin{prp}\label{fff}If a group $G$ satisfies both {\rm Min-n} and {\rm Max-n} conditions and for every $g\in G$ the class $[e]_g$ is a subgroup of $G$, then $G$ is nilpotent.
\end{prp}
From Proposition \ref{fff} follows, in particulat, that if $G$ is a finite group such that for every $g\in G$ the class $[e]_g$ is a subgroup of $G$, then $G$ is nilpotent. 

In the conclusion of this Section we mention that the following commutator identities hold in every group.
\begin{align}
\notag[xy,z]&=[x,z]^y~[y,z]\\
\notag[x,yz]&=[x,z]~[y,z]^x\\
\notag[x,y]^z&=[x^z,y^z]
\end{align}
In particular, from this equalities follows that for every element $x$ which belongs to the subgroup of $G$ generated by the elements $x_1,\dots,x_n$ and for every element $y$ which  belongs to the subgroup of $G$ generated by the elements $y_1,\dots,y_m$ there exist some elements $a_1$, $\dots$, $a_n$, $b_1$, $\dots$, $b_n$, $c_1$, $\dots$, $c_n$, $d_1$, $\dots$, $d_n$ in $G$ such that
\begin{align}
\label{com1}[x,y]&=\prod_{i=1}^n[x,y_i]^{a_i}[x,y_i^{-1}]^{b_i},\\
\label{com2}[x,y]&=\prod_{i=1}^n[x_i,y]^{c_i}[x_i^{-1},y]^{d_i}.
\end{align}  

\section{Structure of the lower central series}\label{swid}

\begin{ttt}\label{wfin}Let $G$ be a group with generators $x_1,\dots,x_n$ such that the twisted conjugacy class $[e]_{g}$ is a subgroup of $G$ for every $g\in G$. Then 
\begin{enumerate}
\item $\gamma_2(G)=[e]_{x_1}[e]_{x_2}\dots[e]_{x_{n-1}}$,
\item if $\gamma_k(G)=[e]_{h_1}[e]_{h_2}\dots[e]_{h_s}$ and $p_{ij}=[x_i,h_j]$, then $\gamma_{k+1}(G)=\prod_{i,j}[e]_{p_{ij}}$.
\end{enumerate}
\end{ttt}
\textbf{Proof.} By Lemma \ref{prim} and Lemma \ref{vnesh} for every elements $y,z\in G$ the commutator $[y,z]$ belongs to $[e]_{y}\cap[e]_z<[e]_y\subseteq[e]_{x_1}[e]_{x_2}\dots[e]_{x_n}[e]_{x_1^{-1}}[e]_{x_2^{-1}}\dots[e]_{x_n^{-1}}$. Since in every group the equality
$$[x,y^{-1}]=\left([x,y]^{y^{-1}}\right)^{-1}$$ 
holds, the subgroups $[e]_{x_i}$ and $[e]_{x_i^{-1}}$ coincide, and hence $[y,z]\in[e]_{x_1}\dots[e]_{x_n}$. From equality (\ref{com2}) and equality $[e]_{x_i}=[e]_{x_i^{-1}}$ follows that $[e]_{x_n}\leq[e]_{x_1}\dots[e]_{x_{n-1}}$, therefore every commutator $[x,y]$ belongs to $[e]_{x_1}\dots[e]_{x_{n-1}}$. Since by Lemma \ref{nsub} for every $g,h\in G$ the classes $[e]_g$, $[e]_h$ are normal subgroups of $G$, we have $[e]_g[e]_h=[e]_h[e]_g$ and therefore every product of commutators belongs to $[e]_{x_1}\dots[e]_{x_{n-1}}$, what proves the first item of the theorem.

Let $\gamma_k(G)=[e]_{h_1}\dots[e]_{h_s}$, then for an arbitrary element $x\in\gamma_k(G)$ there exist some elements $y_1,\dots,y_s$ in $G$ such that $x=[y_1,h_1]\dots[y_s,h_s]$. Using equality (\ref{com2}) for every $i=1,\dots,s$ we have
$$
[y_i,h_i]=[x_1,h_i]^{a_{1i}}[x_2,h_i]^{a_{2i}}\dots[x_n,h_i]^{a_{ni}}[x_1^{-1},h_i]^{b_{1i}}[x_2^{-1},h_i]^{b_{2i}}\dots[x_n^{-1},h_i]^{b_{ni}}
$$
for some elements $a_{1i}, \dots, a_{ni}, b_{1i}, \dots, b_{ni}$ from $G$. If we denote by $p_{ij}=[x_i,h_j]$, $q_{ij}=[x_i^{-1},h_j]$, then we conclude that every element $x$ from $\gamma_k(G)$ can be written in the form
\begin{equation}\label{eq1}
x=\prod_{\substack{
i=1,\dots,n \\
j=1,\dots,s
}}p_{ij}^{a_{ij}}q_{ij}^{b_{ij}}
\end{equation}
for some elements $a_{ij}, b_{ij}~(i=1,\dots,n,j=1,\dots,s)$ of $G$.
For an arbitrary element $t$ from $G$ and element $x$ from $\gamma_k(G)$ using equalities (\ref{com1}) and (\ref{eq1}) we have
$$[t,x]=\bigl[t,\prod_{i,j}p_{ij}^{a_{ij}}q_{ij}^{b_{ij}}\bigr]=\prod_{i,j}\bigl[t,p_{ij}^{a_{ij}}\bigr]^{c_{ij}}\bigl[t,q_{ij}^{b_{ij}}\bigr]^{d_{ij}}=\prod_{i,j}\bigl[t^{a_{ij}^{-1}},p_{ij}\bigr]^{a_{ij}c_{ij}}\bigl[t^{b_{ij}^{-1}},q_{ij}\bigr]^{b_{ij}d_{ij}}$$
for some elements $c_{ij}, d_{ij}~(i=1,\dots,n,j=1,\dots,s)$ from $G$. Therefore since subgroups $[e]_{p_{ij}}$ and $[e]_{q_{ij}}$ are normal in $G$, the element $[t,x]$ belongs to $\prod_{i,j}[e]_{p_{ij}}[e]_{q_{ij}}$ for every $t\in G$, $x\in \gamma_k(G)$. Since in every group the equality 
$$[x^{-1},y]=\left([x,y]^{x^{-1}}\right)^{-1}$$
holds, we have
$$[z,q_{ij}]=[z,[x_i^{-1},h_j]]=\left[z,\left([x_i,h_j]^{x_i^{-1}}\right)^{-1}\right]=\left[z^{x_i},p_{ij}^{-1}\right]^{x_{i}^{-1}}=\left(\left[z^{x_i},p_{ij}\right]^{p_{ij}^{-1}x_i^{-1}}\right)^{-1},$$
therefore $[e]_{q_{ij}}=[e]_{p_{ij}}$ and for every elements $t\in G$, $x\in\gamma_k(G)$ the commutator $[t,x]$ belongs to $\prod_{i,j}[e]_{p_{ij}}$, and  hence $\gamma_{k+1}(G)=\prod_{i,j}[e]_{p_{ij}}$.\hfill$\square$
\begin{cor}\label{lwid}
Let $G$ be a group with generators $x_1,\dots,x_n$ such that the class $[e]_{g}$ is a subgroup of $G$ for every $g\in G$. Then 
\begin{enumerate}
\item ${\rm wid}(\gamma_2(G))\leq n-1$,
\item ${\rm wid}(\gamma_{k}(G))\leq n^{k-2}(n-1)/2$ for $k\geq 3$.
\end{enumerate}
\end{cor}
\textbf{Proof.} The first item of Theorem \ref{wfin} implies that ${\rm wid}(\gamma_2(G))\leq n-1$. From the second item of Theorem \ref{wfin} we have $\gamma_3(G)=\prod_{i,j}[e]_{[x_i,x_j]}$, where $i$ varies from $1$ to $n$ and $j$ varies from $1$ to $n-1$.
Since $[e]_{[x_i,x_i]}=e$ and $[e]_{[x_i,x_j]}=[e]_{[x_i,x_j]^{-1}}=[e]_{[x_j,x_i]}$ we can rewrite $\gamma_3(G)$ in the form
$$\gamma_3(G)=\prod_{\substack{
i=1,\dots,n \\
j=1,\dots,n-1\\
i>j
}}[e]_{[x_i,x_j]}.$$
Thereforel ${\rm wid}(\gamma_3(G))$ is less then or equal to the number of members in this product which is equal to $n(n-1)/2$. From the second item of Theorem \ref{wfin} it follows that ${\rm wid}(\gamma_{k+1}(G))\leq n\cdot{\rm wid}(\gamma_k(G))$, i.~e. ${\rm wid}(\gamma_{k}(G))\leq n^{k-2}(n-1)/2$ for $k\geq 3$.\hfill$\square$
\begin{rmk}{\rm If $G$ is generated by two elements, then ${\rm wid}(\gamma_2(G))\leq2^1-2^0=1$, i.~e. every element of $\gamma_2(G)$ is a commutator. Moreover ${\rm wid}(\gamma_3(G))\leq2(2-1)/2=1$, i.~e. every element from $\gamma_3(G)$ has the form $[[x,y],z]$. Since two generated group $G$ such that $[e]_g\leq G$ for all $g\in G$ does not have to be abelian or $2$-step nilpotent (see Section \ref{exa}), the estimations ${\rm wid}(\gamma_2(G))\leq n-1$ and ${\rm wid}(\gamma_3(G))\leq n(n-1)/2$ cannot be improved.}
\end{rmk}
\begin{problem} {\rm Construct the sharp estimation of  ${\rm wid}(\gamma_k(G))$ for $k\geq4$.}
\end{problem}

If the group $G$ from Corollary \ref{lwid} is metabelian, then for every elements $y_1,\dots,y_m$ $(m\geq3)$ from $G$ and every permutation  $\sigma \in {\rm Sym}(\{3,\dots,m\})$ the commutator equality $[y_1,y_2,y_3,\dots,y_n]=[y_1,y_2,y_{\sigma(3)},\dotsm,y_{\sigma(n)}]$ holds \cite[\S 1.1]{MikPas}. So, in the case of metabelian groups we can strongly improve the estimation for ${\rm wid}(\gamma_{k+2}(G))$ $(k\geq1)$ in the form
$${\rm wid}(\gamma_{k+2}(G))\leq \frac{n(n-1)}{2}\begin{pmatrix}n+k-2\\k-1\end{pmatrix},$$
where $\begin{pmatrix}n+k-2\\k-1\end{pmatrix}$ is the number of ways to distribute $k-1$ indistinguishable balls to the $n$ distinguishable baskets: $\gamma_{k+2}(G)=\prod_{i,j,i_3,\dots,i_{k+1}}[e]_{[x_{i},x_{j},x_{i_3},\dots,x_{i_{k+1}}]}$ and $[e]_{[x_{i},x_{j},x_{i_3},\dots,x_{i_{k+1}}]}$ depents only on $x_{i}$, $x_{j}$ for $i=1,\dots,n$, $j=1,\dots,n-1$ for $i>j$ ($n(n-1)/2$ variants) and the number of each $x_i~(i=1,\dots,n)$ ($n$ distinguishable baskets) in the set $\{x_{i_3},\dots,x_{i_{k+1}}\}$ ($k-1$ indistinguishable balls).

Another improvement of the estimation of the width of the members of lower central series can be obtained in the case when $\gamma_k(G)=\gamma_{k+1}(G)$ for some $k$.
\begin{prp}\label{structure}
Let $G$ be a  group with $n$ generators such that  $[e]_{g}$ is a subgroup of $G$ for every $g\in G$. Let $k\geq2$ satisfies $\gamma_k(G)=\gamma_{k+1}(G)$. Then for some number $N\leq n^{k-1}(n-1)$ there exist $N$ sequences $\{g_{1i}\}_i$, $\{g_{2i}\}_i$, $\dots$,  $\{g_{Ni}\}_i$ such that $g_{11}, g_{21},\dots,g_{N1}\in\gamma_k(G)$ and for all $j=1,2,\dots$ the following equality holds.
$$\gamma_k(G)=[e]_{[g_{11},\dots,g_{1j}]}[e]_{[g_{21},\dots,g_{2j}]}\dots [e]_{[g_{N1},\dots,g_{Nj}]}$$
\end{prp}
\textbf{Proof.} By Theorem \ref{wfin} there exist some elements $h_1,\dots,h_s$ $(s\leq n^{k-2}(n-1))$ from $\gamma_{k-1}(G)$ such that $\gamma_k(G)=[e]_{h_1}\dots[e]_{h_s}$. If we denote by $p_{ij}=[x_i,h_j]$, then by Theorem \ref{wfin}(2) we have $\gamma_{k+1}(G)=\prod_{i,j}[e]_{p_{ij}}$. Let $N=ns\leq n^{k-1}(n-1)$ and $\{g_{11},g_{21},\dots,g_{N1}\}=\{p_{ij}~|~1\leq i\leq n, 1\leq j\leq s\}$ (we just renamed the elements $p_{ij}$ in the way which we need in the formulation of the lemma). All the element $g_{11},g_{21},\dots,g_{N1}$ belong to $\gamma_k(G)$.

Since $\gamma_k(G)=\gamma_{k+1}(G)$,  the element $g_{11}$ belongs to  $\gamma_{k+1}(G)=\prod_{i=1}^N[e]_{g_{i1}}$. Therefore there exist some elements $(g_{12}:=)q_{1},q_{2},\dots,q_{N}$ such that $g_{11}=\prod_{i=1}^N[g_{i1},q_{i}]$.
Then by Lemma \ref{prim} we have
$$[e]_{g_{11}}<\prod_{i=1}^{N}[e]_{[g_{i1},q_{i}]}<[e]_{[g_{11},q_{1}]}\prod_{i=2}^N[e]_{g_{i1}}=[e]_{[g_{11},g_{12}]}\prod_{i=2}^N[e]_{g_{i1}}.$$
Therefore $\gamma_k(G)=\gamma_{k+1}(G)=[e]_{[g_{11},g_{12}]}\prod_{i=2}^N[e]_{g_{i1}}$. Now taking $g_{21}$ instead of $g_{11}$ and repeating the reasoning above we conclude that $\gamma_k(G)=[e]_{[g_{11},g_{12}]}[e]_{[g_{21},g_{22}]}\prod_{i=3}^N[e]_{g_{i1}}$. Repeating this idea for all $g_{i1}$ we have
$\gamma_k(G)=\prod_{i=1}^N[e]_{[g_{i1},g_{i2}]}$.
Repeating this proccess $k$ time we conclude that 
$\gamma_k(G)=\prod_{i=1}^N[e]_{[g_{i1},g_{i2},\dots,g_{ik}]}$
 for some elements $g_{ij}(1\leq i\leq N, 1\leq j\leq k)$. This fact finishes the proof.\hfill$\square$

\begin{cor}\label{bounded}Let $G$ be a  group with $n$ generators such that the class $[e]_{g}$ is a subgroup of $G$ for every $g\in G$ and let $k$ satisifies $\gamma_k(G)=\gamma_{k+1}(G)$. Then ${\rm wid}(\gamma_s(G))\leq n^{k-1}(n-1)$ for all $s$.
\end{cor}
\textbf{Proof.} By Corollary \ref{lwid} for $s\leq k+1$ we have ${\rm wid}(\gamma_s(G))\leq n^{s-1}-n^{s-2}\leq n^{k}-n^{k-1}$ and we only need to estimate the width ${\rm wid}(\gamma_s(G))$ for $s>k+1$. From Proposition~\ref{structure} it follows that ${\rm wid}(\gamma_s(G))\leq N\leq n^{k-1}(n-1)$ for $s>k+1$.\hfill$\square$

\begin{rmk}{\rm In Corollary \ref{bounded} we suppose that $\gamma_k(G)=\gamma_{k+1}(G)$ but when we write ${\rm wid}(\gamma_{k+1}(G))$ we also specify the generating set of $\gamma_k(G)=\gamma_{k+1}(G)$ and the equality $\gamma_k(G)=\gamma_{k+1}(G)$ does not imply the equality ${\rm wid}(\gamma_k(G))={\rm wid}(\gamma_{k+1}(G))$.} 
\end{rmk}

\section{Structure of the group}\label{sstruc}
At first we prove one simple general necessary condition which must be satisfied by an infinite group $G$ where the class $[e]_g$ is a subgroup of $G$ for every $g$.
\begin{prp}\label{conjjjj}Let $G$ be an infinite group such that $[e]_g$ is a subgroup of $G$ for every $g$. Then $G$ has infinitely many conjugacy classes.
\end{prp}
\textbf{Proof.} Suppose that $g$ has $s<\infty$ conjugacy classes, and let $g_1,g_2,\dots$ be elements of $G$. Consider the sequence of elements
$$h_1=g_1, h_2=[g_1,g_2], h_3=[g_1,g_2,g_3],\dots,h_k=[g_1,\dots,g_k],\dots$$
Since there are only $s$ conjugacy classes in the group $G$, there exist two numbers $m<n\leq s+1$ such that $h_m$ and $h_n$ are conjugated, i.~e. $h_m=h_n^{a}$ for some $a\in G$. Therefore 
$$[e]_{h_m}=[e]_{h_n^a}=[e]_{h_n}^a=[e]_{h_n}.$$
The element $h_n=[g_1,\dots,g_n]=[h_m,g_{m+1},\dots,g_n]$ belongs to $[e]_{h_m}$. If $h_n\neq e$, then by Lemma \ref{easy1} it does not belong to $[e]_{h_n}=[e]_{h_m}$. This contradiction says that $h_n=e$. Therefore for $r\geq s+1$ we have $h_r=[h_n,g_{n+1},\dots,g_r]=e$ and since $g_1,g_2,\dots$ is an arbitrary sequence of elements from $G$, then $G$ is nilpotent of nilpotency class less then or equal to $s$. Since a nilpotent group has finitely many conjugacy classes if and only if it is finite, $G$ must be a finite group, what contradicts the condition of the proposition.\hfill$\square$

If $G$ is a groups such that the class $[e]_g$ is a subgroup of $G$ for every $g\in G$, then for elements $g_1, g_2,\dots$ from $G$ if $[e]_{[g_1,\dots,g_k]}$ is not trivial, then 
$$[e]_{g_1}>[e]_{[g_1,g_2]}>[e]_{[g_1,g_2,g_3]}>\dots$$
 is a strictly descending chain of normal subgroups. Then it is reasonable to ask about the structure of $\cap_{k}[e]_{[g_1,\dots,g_k]}$.
\begin{problem}\label{intrat}{\rm Let $G$ be a group such that the class $[e]_g$ is a subgroup of $G$ for every $g$. Is it true that $\cap_{k}[e]_{[g_1,\dots,g_k]}=e$ for every elements $g_1, g_2,\dots$ from $G$?}
\end{problem} 
A positive solution of Problem \ref{intrat} would emply the fact that the lower central series of a group is strictly descending.
\begin{ttt}\label{interest}Let $G$ be a finitely generated group such that the class $[e]_g$ is a subgroup of $G$ for every $g$. If for every sequence of elements $g_1,g_2,\dots$ the intersection  $\cap_{k}[e]_{[g_1,\dots,g_k]}$ is trivial, then either $\gamma_k(G)=e$ or $\gamma_k(G)\neq \gamma_{k+1}(G)$.
\end{ttt}
\textbf{Proof.} Suppose that $\gamma_k(G)=\gamma_{k+1}(G)$, from Proposition \ref{structure} we have 
$$\gamma_k(G)=[e]_{[g_{11},\dots,g_{1j}]}[e]_{[g_{21},\dots,g_{2j}]}\dots [e]_{[g_{N1},\dots,g_{Nj}]},$$
for an infinite set of elements elements $g_{ij}$ and every $j>1$. Suppose that in $\gamma_k(G)$ there exists an element $a\neq e$. Since $\cap_j[e]_{[g_{11},\dots,g_{1j}]}=e$, for $j_1$ big enough the element $a$ does not belong to $[e]_{[g_{11},\dots,g_{1j_1}]}$. Consider the natural homomorphism 
$$\varphi_1:G\to H=G/[e]_{[g_{11},\dots,g_{1j_1}]}.$$
 Since $a$ does not belong to $[e]_{[g_{11},\dots,g_{1j_1}]}$, then by Lemma \ref{quu}
 $$e\neq\varphi_1(a)\in\varphi_1(\gamma_k(G))\leq [e]_{[\varphi_1(g_{21}),\dots,\varphi_1(g_{2j})]}\dots [e]_{[\varphi_1(g_{N1}),\dots,\varphi_1(g_{Nj})]}$$
for $j>j_1$.  Using the same argument we see that there exists a homomorphism 
 $$\varphi_2:H\to H/[e]_{[\varphi_1(g_{21}),\dots,\varphi_1(g_{2j_2})]}$$
such that $e\neq\varphi_2\varphi_1(a)\in\varphi_2\varphi_1(\gamma_k(G))$ and
$$\varphi_2\varphi_1(\gamma_k(G))\leq [e]_{[\varphi_2\varphi_1(g_{31}),\dots,\varphi_2\varphi_1(g_{3j})]}\dots [\varphi_2\varphi_1(e)]_{[\varphi_2\varphi_1(g_{N1}),\dots,\varphi_2\varphi_1(g_{Nj})]}$$
for $j>j_2>j_1$. Repeating this process $N-1$ times we conclude that there exist a homomorphism $\varphi=\varphi_{N-1}\dots\varphi_1$ from $G$ to $G/H$ for some subgroup $H$ such that $e\neq \varphi(a)\in\varphi(\gamma_k(G))\leq [e]_{[\varphi(g_{N1}),\dots,\varphi(g_{Nj})]}$ for $j$ big enough. The equality $\cap_k[e]_{[g_{N1},\dots,g_{Nk}]}=e$ implies $\varphi(a)=e$ what contradicts the choice of $a$, i.~e. in $\gamma_k(G)$ there are no nontrivial elements.\hfill $\square$

\begin{cor}\label{minn} Let $G$ be a finitely generated group such that $[e]_g$ is a subgroup of $G$ for every $g\in G$. If $G$ satisfies {\rm Min-n} condition, then $G$ is nilpotent.
\end{cor}
\textbf{Proof.} From Min-n condition follows that for every sequence of elements $\{g_i\}_i$ the intersection  $\cap_{k}[e]_{[g_1,\dots,g_k]}$ is trivial. If $\gamma_k(G)=\gamma_{k+1}(G)$ for some $k$, then by Theorem \ref{structure} we have
$\gamma_k(G)=e$, i.~e. $G$ is nilpotent.
If $\gamma_{k+1}(G)\neq\gamma_k(G)$ for all $k$, then $\gamma_1(G)>\gamma_2(G)>\dots$ is a strictly descending chain of normal subgroups, therefore by Min-n condition $\gamma_k(G)=e$ for some $k$.
\hfill$\square$

Corollary \ref{minn} generalizes the result from \cite[Theorem 3]{BarNasNes} which we recalled in Proposition \ref{fff}: we change Max-n condition by a weaker condition which says that $G$ is a finitely generated group.

\section{Finite abelian-by-cyclic groups}
Let a group $H$ acts on a group $A$, i.~e. there exists a homomorphism from $H$ to ${\rm Aut}~A$ which maps an element $h$ from $H$ to the automorphism $\hat{h}:a\mapsto a^h$ of the group $A$. The map $d:H\to A$ is called \textit{a derivation} if for every elements $h_1,h_2\in H$ the equality 
\begin{equation}\label{derr}
d(h_1h_2)=d(h_1)d(h_2)^{h_1}
\end{equation}
 holds. The set of all derivations $H\to A$ is denoted by the symbol $Der(H,A)$.

Let $H=\mathbb{Z}_n=\langle t\rangle$ be a cyclic group of order $n$ and $A$ be an abelian group (with the operation written additively). Denote by $\psi$ the automorphism $a\mapsto a^t$. If $d$ is a derivation from $H$ to $A$ such that $d(t)=a$, then using equality (\ref{derr}) we conclude that
\begin{equation}\label{powers}
d(t^k)=(id+\psi+\psi^2+\dots+\psi^{k-1})a
\end{equation}
for all $k=1,\dots,n$. So, in order to define a derivation from $\mathbb{Z}_n=\langle t\rangle$ to $A$ it is enough to define only the valude $d(t)=a$. The map $d(t)=a$ can be extended to the derivation $\mathbb{Z}_n=\langle t\rangle\to A$ if an only if
\begin{equation}\label{gennecsof}
0=d(e)=d(t^n)=(id+\psi+\psi^2+\dots+\psi^{n-1})a.
\end{equation}
In particular, for every $b\in A$ the map $d(t)=(id-\psi)b$ can be extended to the derivation from $Der(\mathbb{Z}_n,A)$. Such derivations are called \textit{principal derivations} \cite[Chapters II, III]{Bro}.
 
Let $p>2$ be a prime, $G$ be a finite $p$-group, $A$ be a normal abelian subgroup of $G$ and $H=G/A$. For $g\in G$ denote by $\overline{g}=gA$ the image of $g$ under the natural homomorphism $G\to H$. If $d\in Der(H,A)$, then the map
$$\varphi_d:g\mapsto d(\overline{g})g$$
is an automorphism of $G$ (see \cite[Satz 2]{Gas} for details).

Following \cite{Khu} for a finite abelian $p$-group $A$ denote by $\Omega_k(A)=\{a\in A~|~a^{p^k}=e\}$.
\begin{lem}\label{derivations} Let $p>2$ be a prime, $A$ be a finite abelian $p$-group and $G$ be an extension of $A$ by $\mathbb{Z}_{p^n}=\langle t\rangle$. If the class $[e]_{\varphi}$ is a subgroup of $G$ for every $\varphi\in {\rm Aut}~G$, then for every $a\in\Omega_1(A)$ the map $d(t)=a$ can be extended to the derivation $d\in Der(\mathbb{Z}_{p^n},A)$.
\end{lem}
\textbf{Proof.} Denote by $\psi$ the automorphism $\psi:a\mapsto a^t$ of the group $A$. Since $\psi^{p^n}=id$, the order of $\psi$ is a power of $p$. Using the additive form for the operation in abelian group $A$, by equality (\ref{gennecsof}) the map $d(t)=a$ can be extended to derivation $d:\mathbb{Z}_{p^n}\to A$ if an only if
\begin{equation}\label{necsof}
(id+\psi+\dots+\psi^{p^n-1})a=0.
\end{equation}
So, we need to prove that for every element $a\in \Omega_1(A)$ the equality (\ref{necsof}) holds.

Let $c$ be an arbitrary element of the group $A$ and $b=(id-\psi)c$. The map $d(t)=b$ can be extended to the derivation $d:\mathbb{Z}_{p^n}\to A$ and by the  condition of the lemma, for the automorphism $\varphi:g\mapsto d(\overline{g})g$ the class $[e]_{\varphi}$ is a subgroup of $G$. This subgroup consists of the elements $g\varphi(g)^{-1}=d(\overline{g})^{-1}$. Since $[e]_{\varphi}$ is a subgroup, we have $[e]_{\varphi}=[e]_{\varphi}^{-1}=d(\mathbb{Z}_{p^n})$ and by equality (\ref{powers}) this subgroup consists of elements
$$(id+\psi+\dots+\psi^k)b$$
for all $k=1,\dots,p^{n}$. Since by Lemma \ref{nsub} the class $[e]_{\varphi}$ is a normal subgroup of $G$, the map $\psi$ maps $[e]_{\varphi}$ to itself. Therefore since $b$ belongs to $[e]_{\varphi}$, the elements $\psi^k(b)$ belong to $[e]_{\varphi}$ for all $k=1,\dots,p^n$. Since the order of $[e]_{\varphi}=d(\mathbb{Z}_{p^n})$ is less then or equal to $p^n$ and all elements $\psi^k(b)~(k=1,\dots,p^n)$ are not trivial, all the elements $\psi^k(b)~(k=1,\dots,p^n)$ cannot be distinct. Therefore $\psi^r(b)=\psi^s(b)$ for some $0\leq r\neq s\leq p^n-1$, i.~e. $\psi^{k}(b)=b$ for some $0< k=r-s\leq p^n-1$. Since the order of $\psi$ is a power of $p$, we can assume that $k$ is a power of $p$, i.~e.  for some $0\leq k\leq n-1$ we have $\psi^{p^k}(b)=b$. Acting on this equality by $\psi^{p^k(p^{n-k-1}-1)}$ we have $\psi^{p^{n-1}}(b)=b$ for all $b=(id-\psi)c$. Since $c$ is an arbitrary element of $A$ we have the equality $\psi^{p^{n-1}}(id-\psi)=(id-\psi)$ which emplies the following equality
$$\psi^{p^{n-1}+1}=\psi^{p^{n-1}}+\psi-id.$$
Using simple induction on the pair $(r,s)$ in lexicographic order, this equality  implies the equality
$$\psi^{rp^{n-1}+s}=r(\psi^{p^{n-1}}-id)+\psi^s$$
for all $r=0,\dots,p-1$ and $s=0,\dots,p^{n-1}-1$. Therefore we have
\begin{align}
\notag id+\psi+\dots+\psi^{p^n-1}&=(id+\psi+\dots+\psi^{p^{n-1}-1})\sum_{r=0}^{p-1}\psi^{rp^{n-1}}\\
\notag&=(id+\psi+\dots+\psi^{p^{n-1}-1})\sum_{r=0}^{p-1}r(\psi^{p^{n-1}}-id)\\
\notag&=(id+\psi+\dots+\psi^{p^{n-1}-1})(\psi^{p^{n-1}}-id)\frac{p(p-1)}{2}
\end{align}
This equality emplies that $(id+\psi+\dots+\psi^{p^n-1})a=0$ for every element $a$ of order $p$ (i.~e. $a\in\Omega_1(A)$), so, the map $d(t)=a$ can be extended to the derivation $d:\mathbb{Z}_{p^n}\to A$.\hfill$\square$
\begin{lem}\label{om1} Let $p>2$ be a prime, $A$ be a finite abelian $p$-group and $G$ be an extension of $A$ by $\mathbb{Z}_{p^n}=\langle t\rangle$. If the class $[e]_{\varphi}$ is a subgroup of $G$ for every $\varphi\in {\rm Aut}~G$, then $a^t=a$ for every $a\in\Omega_1(A)$.
\end{lem}
\textbf{Proof.} Denote by $\psi$ the automorphism $\psi:a\mapsto a^t$ of the group $A$. Since $\psi^{p^n}=id$, the order of $\psi$ is a power of $p$. Let $a\in \Omega_1(A)$ be an arbitrary element. Then by Lemma \ref{derivations} the map $d(t)=a$ defines a derivation. For the automorphism $\varphi:g\mapsto d(\overline{g})g$ the class $[e]_{\varphi}$ is a subgroup of $G$ which coincides with $d(\mathbb{Z}_{p^n})$ and consists of the elements
\begin{equation}\label{elementy}
d(t^k)=(id+\psi+\dots+\psi^{k-1})a
\end{equation}
for $k=1,\dots,p^n$. Similarly to Lemma \ref{derivations} we have $\psi^{p^m}(a)=a$ for some $0\leq m\leq n-1$. If $m=0$, then there is nothing to prove. Suppose that $\psi(a)\neq a$ and let $m\neq0$ be a minimal number such that $\psi^{p^m}(a)=a$.

Suppose that $(id+\psi+\dots+\psi^{p^m-1})a=0$. In this case using the equality (\ref{elementy}) the group $d(\mathbb{Z}_{p^n})$ consists of the elements
$$(id+\psi+\dots+\psi^{k-1})a$$
for $k=1,\dots,p^m$, i.~e. it contains only $p^m$ elements. This group contains $a$ and it is fixed under the action of $\psi$. Therefore $\psi^r(a)$ belongs to $d(\mathbb{Z}_{p^n})$ for all $r=1,\dots,p^m$. Since $d(\mathbb{Z}_{p^n})$ contains only $p^m$ elements and all the elements $\psi^r(a)$ are not trivial, there exists some numbers $0\leq r\neq s<p^m$ such that $\psi^r(a)=\psi^s(a)$, i.~e. $\psi^k(a)=a$ for some $0<k=r-s<p^m$. Since the order of $\psi$ is a power of $p$, we can assume that $k$ is a power of $p$. It contradicts the fact that $m$ is a minimal number such that $\psi^{p^m}(a)=a$. Therefore $(id+\psi+\dots+\psi^{p^m-1})a\neq0$.

Let $b=(id+\psi+\dots+\psi^{p^m-1})a\neq0$ and $H$ be a group generated by $b$. From the equality $\psi^{p^m}(a)=a$ follows that
$$(id+\psi+\dots+\psi^{rp^{m}-1})a=r(id+\psi+\dots+\psi^{p^m-1})a,$$
for all $r=1,\dots,p^{n-m}$. Therefore the group $H$ is a subgroup of $[e]_{\varphi}$. Since the element $b=(id+\psi+\dots+\psi^{p^m-1})a$ is fixed under $\psi$, the group $H$ is a normal subgroup of $G$ which  does not contain $a$. The quotient group $[e]_{\varphi}/H$ consists of the elements
$$(id+\psi+\dots+\psi^{k-1})a+H$$
for $k=1,\dots,p^m$. This subgroup contains  $a+H$ and $\psi([e]_{\varphi}/H)=[e]_{\varphi}/H$. Therefore the elements $\psi^r(a)+H$ belong to $[e]_{\varphi}/H$ for all $r=0,\dots,p^{m}-1$. Since the group $[e]_{\varphi}/H$ contains less then or equal to $p^m$ elements and all the elements $\psi^r(a)+H$ are not trivial in $[e]_{\varphi}/H$, there exist some numbers $0\leq r\neq s<p^{m}$ such that $\psi^r(a)+H=\psi^s(a)+H$. So, $\psi^k(a)+H=a+H$  for some $0<k=r-s\leq p^m-1$. Since the order of $\psi$ is a power of $p$, we can assume that $k=p^l$ for some $0\leq l\leq m-1$. So, for some number $r$ we have $\psi^{p^l}(a)=a+r(id+\psi+\dots+\psi^{p^m-1})a$. Since $(id+\psi+\dots+\psi^{p^m-1})a$ is fixed under the action of $\psi$, then for all $s$ we have $\psi^{sp^l}(a)=a+sr(id+\psi+\dots+\psi^{p^m-1})a$.
Therefore
\begin{align}
\notag \psi^{p^l}(a)&=a+r(id+\psi+\dots+\psi^{p^m-1})a\\
\notag &=a+r(id+\psi+\dots+\psi^{p^l-1})\left(\sum_{s=0}^{p^{m-l}-1}\psi^{sp^{l}}\right)a\\
\notag &=a+r(id+\psi+\dots+\psi^{p^l-1})\sum_{s=0}^{p^{m-l}-1}(a+sr(id+\psi+\dots+\psi^{p^m-1})a)\\
\notag&=a+r(id+\psi+\dots+\psi^{p^l-1})p^{m-l}a\\
\notag&+r(id+\psi+\dots+\psi^{p^l-1})(id+\psi+\dots+\psi^{p^m-1})\frac{p^{m-l}(p^{m-l}-1)}{2}a=a,
\end{align}
what contradicts the fact that $m$ is a minimal number such that $\psi^{p^m}(a)=a$. So, $\psi(a)$ must be equal to $a$.
\hfill$\square$
\begin{cor}\label{omk} Let $p>2$ be a prime, $A$ be a finite abelian $p$-group and $G$ be an extension of $A$ by $\mathbb{Z}_{p^n}=\langle t\rangle$. If the class $[e]_{\varphi}$ is a subgroup of $G$ for every $\varphi\in {\rm Aut}~G$, then $a^t=a\theta(a)$, where $\theta$ is a homomorphism from $A$ to $\Omega_n(A)$ such that $\theta(\Omega_k(A))\leq\Omega_{k-1}(A)$ for all $k\leq n$.
\end{cor}
\textbf{Proof.} Denote by $\psi$ the automorphism $\psi:a\mapsto a^t$ of the group $A$. Let $b$ be an arbitrary element of the group $A$ and  $a=(id-\psi)b$. The map $d(t)=a$ defines a derivation from $Der(\mathbb{Z}_{p^n},A)$ and for the automorphism $\varphi:g\mapsto d(\overline{g})g$ the class $[e]_{\varphi}=d(\mathbb{Z}_{p^n})$ is a subgroup of $G$ which contains $a$. Since the order of $[e]_{\varphi}=d(\mathbb{Z}_{p^n})$  is less then or equal to $p^n$, then the order of $a$ is less then or equal to $p^n$, i.~e. for every $b$ the element $(id-\psi)b$ belongs to $\Omega_n(A)$. So, for every $b\in A$ we have $b^t=\psi(b)=b\theta(b)$, where $\theta$ is the map from $A$ to $\Omega_n(A)$. For arbitrary elements $a,b\in A$ we have 
$\theta(ab)ab=\psi(ab)=\psi(a)\psi(b)=a\theta(a)b\theta(b)=ab\theta(a)\theta(b)$, therefore $\theta$ is a homomorphism. 

If $a$ is an element of $\Omega_k(A)$, then $a^{p^{k-1}}$ is an element from $\Omega_1(A)$, and then by Lemma \ref{om1} we have $\theta(a)^{p^{k-1}}=\theta(a^{p^{k-1}})=e$, so, $\theta(a)$ belongs to $\Omega_{k-1}(A)$.
\hfill$\square$
\begin{ttt}\label{corex} Let $n$ be an odd integer with prime decomposition $n=p_1^{n_1}\dots p_m^{n_m}$ and $A$ be a finite abelian group. If $G$ is an extension of $A$ by $\mathbb{Z}_n$ such that the class $[e]_{\varphi}$ is a subgroup of $G$ for every $\varphi\in {\rm Aut}~G$, then $G$ is nilpotent of nilpotency class at most ${\rm max}\{n_1,\dots,n_m\}+1$.
\end{ttt}
\textbf{Proof.} Since $G$ is a finite group such that the class $[e]_{g}$ is a subgroup of $G$ for every $g\in G$, by Corollary \ref{minn} the group $G$ is nilpotent. By Burnside-Wielandt theorem \cite[Theorem 17.1.4]{Fund} the group $G$ is the direct product of its Sylow subgroups $G_1, \dots, G_k$. We can assume that $k=m$ (otherwise rewrite $n=p_1^{n_1}\dots p_m^{n_m}p_{m+1}^{0}\dots p_k^0$). 

Since $G$ is an extension of the abelian group $A$ by the cyclic group $\mathbb{Z}_n$, the group $G_i$ must be an extension of the abelian $p_i$-group $A_i$ by the cyclic group $\mathbb{Z}_{p_i^{n_i}}=\langle t_i\rangle$. Denote by $\psi_i$ the automorphism $a\mapsto a^{t_i}$ of the group $A_i$. By Corollary \ref{omk} for every $a\in A_i$ we have $\psi_i(a)=a\theta_i(a)$. Using simple induction on $r\geq2$  it is easy to see that 
$$\gamma_r(G_i)\leq \theta_i^{r-1}(A)$$
and by Corollary \ref{omk} we have $\gamma_{n_i+2}(G_i)=e$, i.~e. $G_i$ is a nilpotent group of nilpotency class at most $n_i+1$ and $G$ is nilpotent group of nilpotency class at most ${\rm max}\{n_1,\dots,n_m\}+1$.
\hfill$\square$

\begin{rmk} {\rm In Lemma \ref{derivations}, Lemma \ref{om1}, Corollary \ref{omk} and Theorem \ref{corex} we require that the class $[e]_{\varphi}$ is a subgroup of $G$ for every automorphism $\varphi\in {\rm Aut}~G$. However in proofs we use only inner automorphisms and automorphisms of the form $g\mapsto d(\overline{g})g$ for $d\in Der(\mathbb{Z}_n,A)$. So, in all of this statements it is enough to require that the class $[e]_{\varphi}$ is a subgroup of $G$ only for these specific automorphisms.}
\end{rmk}
If the group $A$ in Theorem \ref{corex} is elementary abelian, then by Corollary \ref{omk} the group $G$ must be abelian. So we have the following problem.
\begin{problem}{\rm Let $p$ be an odd prime, $A$ be a finite abelian group and $G$ be an extension of $A$ by $\mathbb{Z}_p$ such that the class $[e]_{\varphi}$ is a subgroup of $G$ for every $\varphi\in {\rm Aut}~G$. Is it true that $G$ is abelian?}
\end{problem}
\section{Examples and counterexamples}\label{exa}
Conjecture \ref{prob1} says that if $G$ is a group such that the class $[e]_{\varphi}$ is a subgroup of $G$ for every $\varphi\in {\rm Aut}~G$, then $G$ is nilpotent. Example 1 and Example 2 below show that it is not sufficient to consider only inner automorphisms to prove Conjecture \ref{prob1}. Also this examples give answers to some questions from \cite{BarNasNes}.

We thank E.~Khukhro, who noticed that in the group $G(p,n)$ from Example~1 the twisted conjugacy class of the unit element is a subgroup for every inner automorphism and showed us this example. We also thank V.~Bludov, who demonstrates us Example 2 in details.

~\\
\noindent \textbf{{\scshape Example 1.}} Let $n>2$ be a positive integer and $p>2$ be a prime. Denote by $G(p,n)$ the semidirect product of two cyclic groups $\mathbb{Z}_{p^n}\rtimes\mathbb{Z}_{p^{n-1}}$ with generators $x=(1,0)$, $y=(0,1)$ and relation $yxy^{-1}=x^{p+1}$. The order of $G(p,n)$ is $p^{2n-1}$, therefore $G(p,n)$ is nilpotent. Since 
$$\underbrace{[[\dots[[x,y^{-1}],y^{-1}],\dots],y^{-1}]}_{m ~\text{pairs of square brackets}}=x^{p^m},$$
the nilpotency class of $G(p,n)$ is greater then or equal to $n-1$. Using simple induction on $r$ we see that in $G(p,n)$ the following equality holds for all positive integers $a,s$.
\begin{equation}\label{nux}
(x^by)^s=\left(b\left(1+(p+1)+(p+1)^2+\dots+(p+1)^{s-1}\right),s\right).
\end{equation}
If $\varphi$ is inner automorphism of $G$, then 
$$\varphi:~~x\mapsto x^a,~~y\mapsto x^by.$$
From \cite[Corollary 4.4 (1b)]{GolGon} it follows that $a\in \mathbb{Z}_{p^n}^*$ ($a$ is an invertible element of the ring $\mathbb{Z}_{p^n}$) and $b\equiv0~({\rm mod}~p)$.

The twisted cojugacy class $[e]_{\varphi}$ of the unit element $e$ consists of the elements $(x^ry^s)\varphi(x^ry^s)^{-1}=(x^ry^s)((x^ay^b)^r(x^cy)^s)^{-1}$ for $r=0,\dots,p^{n}-1$, $s=0,\dots,p^{n-1}-1$. Using equality (\ref{nux}) and direct calculations we conclude that the class $[e]_{\varphi}$ consists of the elements 
$$\left(r(1-a)-b\left(1+(p+1)+\dots+(p+1)^{s-1}\right),0\right),$$
for $r=0,\dots,p^{n}-1$, $s=0,\dots,p^{n-1}-1$. Since $b\equiv0~({\rm mod}~p)$, we have $b=kp$ and  $b\left(1+(p+1)+\dots+(p+1)^{s-1}\right)=kp\left(1+(p+1)+\dots+(p+1)^{s-1}\right)=k[(p+1)^s-1]$. Since the multiplicative order of $p+1$ in $\mathbb{Z}_{p^n}$ is equal to $p^{n-1}$ \cite [Chapter IV \S 3.5]{Zas}, the set $\{(p+1)^s-1~|~s=0,\dots,p^{n-1}-1\}$ consists of $p^{n-1}$ different elements from  $\mathbb{Z}_{p^n}$ each of which is divisible by $p$. Since in $\mathbb{Z}_{p^n}$ there are only $p^{n-1}$ elements divisible by $p$, the set $\{(p+1)^s-1~|~s=0,\dots,p^{n-1}-1\}$ contains all this elements. Therefore the class $[e]_{\varphi}$ consists of the elements
\begin{equation}\label{7pr}
\left(r(1-a)-kz,0\right),
\end{equation}
for $r=0,\dots,p^{n}-1$, $z$ is divisible by $p$. This set of elements obviously forms a subgroup of $G(p,n)$.

So, $G(p,n)$ is a nilpotent group of nilpotency class greater then or equal to $n-1$ such that the twisted conjugacy class of the unit element is a subgroup of $G(p,n)$ for every inner automorphism of $G(p,n)$. Denote by $p_n$ the $n$-th prime and by $G$ the infinite direct product 
$$G=\prod_{n>2} G(p_n,n).$$ 
Then $G$ is non-nilpotent residually nilpotent group such that the class $[e]_{\varphi}$ is a subgroup of $G$ for every  inner automorphism $\varphi$ of $G$.
\begin{rmk}{\rm If $\varphi$ is not inner automorphism of $G(p,n)$, then $[e]_{\varphi}$ is not neccessarily a subgroup of $G(p,n)$. For example,  by \cite[Corollary 4.4 (1b)]{GolGon} the map 
$$\varphi:~x\to xy^3,~y\mapsto y$$
 induces an automorphism of $G(3,3)$. However, the class $[e]_{\varphi}$ contains $3$ elements $(0,0)$, $(0,3)$, $(9,6)$ which do not form a subgroup. So, group $G$ does not give a counter example to Conjecture \ref{prob1}.}
\end{rmk}

The group $G(p,n)$ for $p>2$, $n>3$ gives a negative answer to the following problems from \cite{BarNasNes}.

~\\
\noindent\textbf{\cite[Problem 2]{BarNasNes}.} Let $G$ be a finite group such that the class $[e]_g$ is a subgroup of $G$ for every $g\in G$. Can we assert that $G$ is a $2$-step nilpotent group?

~\\
\noindent\textbf{\cite[Problem 3]{BarNasNes}.} Let $G$ be an extension of an abelian group by a cyclic (or abelian) group. Can we assert that for some $g\in G$ the class $[e]_g$ is not a subgroup of $G$?

Also using groups $G(p,n)$ from Example 1 we can answer the following question from \cite{BarNasNes}.

~\\
\noindent\textbf{\cite[Problem 7]{BarNasNes}.} Let $G$ be a group such that  $[e]_g$ is a subgroup of $G$ for every $g\in G$. Let $N$ be a normal subgroup of $G$ and denote by $[e]_{g,N}=\{[h,g]~|~h\in N\}$.
\begin{enumerate}
\item Is it true that $[e]_{g,N}$ is a subgroup of $N$?
\item How are $[e]_{g,N}$ and $[e]_g\cap N$ related?
\end{enumerate}
Let $G=G(3,3)=\langle x,y~|~x^{27}=1, y^{9}=1, yxy^{-1}=x^{4}\rangle$ and $N=\langle x^3, y^3\rangle$. Then $[e]_{x,N}$ contains only three elements $e$, $x^{15}$ and $x^{18}$ which do not form a subgroup. It gives negative answer to the first item of the question.

Let $G=G(p,n)=\langle x,y~|~x^{p^n}=1, y^{p^{n-1}}=1, yxy^{-1}=x^{p+1}\rangle$. If $N=G$, then $[e]_{g,N}=[e]_g$, i.~e. $[e]_{g,N}$ and $[e]_g\cap N$ can coincide. If $N=\langle x\rangle=\mathbb{Z}_{p^n}$, then $[e]_{x^{p^k},N}=e$ but from equality (\ref{7pr}) follows that $[e]_{x^{p^k}}\cap N=[e]_{x^{p^k}}=\langle x^{p^{k+1}}\rangle$, therefore  for $[e]_{g,N}=e$ the group $[e]_g\cap N$ can be arbitrary. Therefore there are no essential relations between $[e]_{g,N}$ and $[e]_g\cap N$.

~\\
\noindent \textbf{{\scshape Example 2.}} Let $\mathbb{F}$ be a field, $T$ be a variable and $\mathbb{F}(T)$ be the field of rational function over $\mathbb{F}$ with variable $T$. Denote by $P=\mathbb{F}(T)^+$ the additive group of the field $\mathbb{F}(T)$, by $D=\{d(T)\in \mathbb{F}(T)~|~d(1)=0\}$ and by $H=\{h(T)\in\mathbb{F}(T)~|~h(1)=1\}$. The group $D$ is a subgroup of $P$ and the group $H$ is a subgroup of the multiplicative group $\mathbb{F}(T)^*$ of the field $\mathbb{F}(T)$. Denote by $G$ the group $G=\left\{{\footnotesize \begin{pmatrix}
1&f\\
0&h
\end{pmatrix}}~\biggl|~f\in\mathbb{F}(T), h\in H\right\}$. For elements $u={\footnotesize \begin{pmatrix}
1&f\\
0&h
\end{pmatrix}}$ and $v={\footnotesize \begin{pmatrix}
1&a\\
0&b
\end{pmatrix}}$ we have
\begin{equation}\label{lincom}[u,v]={\begin{pmatrix}
1&a(1-h)+f(b-1)\\
0&1
\end{pmatrix}}.
\end{equation}
If $b\neq1$, then denoting by $h=1$, $f=c(b-1)^{-1}$ we have $$[e]_v=\left\{{ \begin{pmatrix}
1&c\\
0&1
\end{pmatrix}}~\biggl|~c\in\mathbb{F}(T)\right\}\cong P,$$
i.~e. $[e]_v$ is a subgroup of $G$. If $b=1$, then denoting by $h=1-d$ for $d\in D$ we have
$$[e]_v=\left\{{ \begin{pmatrix}
1&ad\\
0&1
\end{pmatrix}}~\biggl|~d\in D\right\}\cong aD,$$
i.~e. $[e]_v$ is again a subgroup of $G$. So the twisted conjugacy class $[e]_g$ is a subgroup of $G$ for every $g\in G$. 

From equality (\ref{lincom}) follows that the center $Z(G)$ is trivial. So, Example 2 gives a negative answer to the following problem from \cite[Problem 18.15]{Kou}.

~\\
\noindent\textbf{\cite[Problem 18.15]{Kou}.} Let $G$ be a group with trivial center. Is it true that there exists an element $g\in G$ such that the class $[e]_g$ is not a subgroup of $G$. 

~\\
Since $Z(G)$ is trivial, the group $G$ is not nilpotent. However, it is residually nilpotent.

\begin{problem}\label{mainp}{\rm Let $G$ be a finitely generated  group such that the class $[e]_{g}$ is a subgroup of $G$ for every $g\in G$. Is it true that $G$ is residually nilpotent?}
\end{problem}

If the group $G$ from Problem \ref{mainp} is residually finite, then by Corollary \ref{minn} it is also residually nilpotent. In particular, if $G$ is finitely generated metabelian or finitely generated linear, then $G$ is residually finite and therefore residually nilpotent.

The positive solution of Problem \ref{mainp} implies the positive solution of Problem \ref{intrat}. However, the positive solution of Problem \ref{intrat} together with Theorem \ref{interest} gives only a hope to have a positive solution of Problem \ref{mainp}.

~\\
\noindent Daciberg Lima Gon\c{c}alves\\
Department of Mathematics-IME, University of S\~{a}o Paulo\\
05508-090, Rua do 
Mat\~{a}o 1010, Butanta-S\~{a}o Paulo-SP, Brazil\\
e-mail: dlgoncal@ime.usp.br\\

\noindent Timur Nasybullov\\
Department of Mathematics, KU Leuven KULAK\\
8500, Etienne Sabbelaan 53, Kortrijk, Belgium\\
e-mail: timur.nasybullov@mail.ru

\end{document}